\documentclass[11pt]{amsart}

\usepackage{amssymb, amsthm, amsmath, gensymb}
\usepackage{graphicx, comment}
\usepackage[small]{caption}
\usepackage{subcaption}
\usepackage{epsfig}
\usepackage{tikz, pgfplots, float}

\usepackage{amsfonts}

\usepackage[utf8]{inputenc}
\usepackage{fullpage}
\usepackage{framed}
\usepackage{multirow}
\usepackage{enumerate}
\usepackage{url}
\usepackage[breaklinks]{hyperref}
\usepackage{cleveref}
\hypersetup{
	colorlinks = true, % colors links instead of ugly boxes
	urlcolor = cyan, % color for external hyperlinks
	linkcolor = teal, % color of internal links
	citecolor = cyan % color of citations
}

\newcommand{\ex}{\mathop{}\!\mathrm{ex}}
\newtheorem{thm}{Theorem}[section]
\newtheorem{definition}[thm]{Definition}
\newtheorem{lem}[thm]{Lemma}
\newtheorem{cor}[thm]{Corollary}
\newtheorem{prop}[thm]{Proposition}
\newtheorem{prob}[thm]{Problem}
\newtheorem{quest}[thm]{Question}
\newtheorem{obs}[thm]{Observation}

\newcommand{\NN}{\mathbb{N}} %  set of natural numbers
\newcommand{\FF}{\boldsymbol{F}}
\newcommand{\HH}{\boldsymbol{H}}
\newcommand{\GG}{\boldsymbol{G}}
\newcommand{\bS}{\boldsymbol{S}}% bold for simplicial complexes

\newcommand{\cF}{{\mathcal F}}
\newcommand{\cE}{{\mathcal E}}

\newcommand{\cD}{\boldsymbol{\mathcal D}}

\newcommand{\cH}{{\mathcal H}}
\newcommand{\cS}{{\mathcal S}}

\newcommand{\cN}{{\mathcal N}}

\newcommand{\cl}{{c\ell}}

\title{Tur\'an problems for simplicial complexes}

\author{Maria Axenovich}
\address{Karlsruhe Institute of Technology, Germany}
\email{maria.aksenovich@kit.edu, liu@mathe.berlin}

\author{D\'aniel Gerbner}
\address{Alfr\'ed R\'enyi Institute of Mathematics, Hungary}
\email{gerbner@renyi.hu, patkos.balazs@renyi.mta.hu}

\author{Dingyuan Liu}

\author{Bal\'azs Patk\'os}

\begin{document}

\maketitle

\begin{abstract}
An abstract simplicial complex $\mathbf{F}$ is a non-uniform hypergraph without isolated vertices, whose edge set is closed under taking subsets. The extremal number $\mathrm{ex}(n,\mathbf{F})$ is defined as the maximum number of edges in an $n$-vertex $\mathbf{F}$-free simplicial complex. Although Tur\'an-type problems for simplicial complexes have long appeared in extremal set theory, a systematic study of $\mathrm{ex}(n,\mathbf{F})$ was initiated only recently by Conlon, Piga, and Sch\"ulke. In contrast to uniform hypergraphs, even the order of magnitude of $\mathrm{ex}(n,\mathbf{F})$ remains unknown for most simplicial complexes.

In this paper, we present a general framework for estimating $\mathrm{ex}(n,\mathbf{F})$ via generalised Tur\'an numbers of associated hypergraphs. Using this approach, we determine the asymptotic behaviour, and in some cases the exact value, of $\mathrm{ex}(n,\mathbf{F})$ for broad classes of simplicial complexes, extending a result of Conlon, Piga, and Sch\"ulke. We also exhibit simplicial complexes whose extremal numbers are not governed by the corresponding generalised Tur\'an numbers, and determine their extremal numbers asymptotically. In addition, we study how the extremal number changes when a new edge is added to the forbidden simplicial complex, and obtain a tight bound for this behaviour.
\end{abstract}

%%%%%%%%%%%%%%%%%%%%%%%%%%%%%%%%%%%%%%%%%%%%%%%%%%%%%%%%%%%%%%%%%%%%%%
%%%%%%%%%%%%%%%%%%%%%%%%%%%%%%%%%%%%%%%%%%%%%%%%%%%%%%%%%%%%%%%%%%%%%%
\section{Introduction}
%%%%%%%%%%%%%%%%%%%%%%%%%%%%%%%%%%%%%%%%%%%%%%%%%%%%%%%%%%%%%%%%%%%%%%
%%%%%%%%%%%%%%%%%%%%%%%%%%%%%%%%%%%%%%%%%%%%%%%%%%%%%%%%%%%%%%%%%%%%%%
An \emph{(abstract) simplicial complex} $\FF=(V,E)$ is a hypergraph, where $V$ is the set of vertices and $E\subseteq2^V$ is the collection of edges that contains all singleton sets and is closed under taking subsets. In other words, $\FF$ is a downward-closed hypergraph without isolated vertices. Note that the empty set is also regarded as an edge of a simplicial complex. The \emph{dimension} of a simplicial complex is defined as its maximum edge size minus $1$. Throughout this paper we consider only simplicial complexes with dimension at least $1$, namely, the largest size of an edge is at least $2$.%\vspace{1mm}

Given a simplicial complex $\FF$, we say that a simplicial complex $\HH$ is \emph{$\FF$-free} if $\HH$ does not contain an isomorphic copy of $\FF$. For $n\in\NN$, the \emph{Tur\'an number} (also called the \emph{extremal number}) $\ex(n,\FF)$ is defined as the maximum number of edges in an $\FF$-free simplicial complex on $n$ vertices.
While the Tur\'an problem for simplicial complexes was considered much earlier, see~\cite[etc.]{Fr83,Sa72,Sh72}, the systematic study was initiated only recently by Conlon, Piga, and Sch\"ulke~\cite{CPS}.%\vspace{1mm}

It is natural to try to relate the extremal number of a simplicial complex $\FF$ to the extremal numbers of some associated hypergraphs. For $s\in[k]$, let $\FF^s$ be the $s$-uniform hypergraph consisting of all the edges in $\FF$ of size $s$. As observed in~\cite{CPS}, for every $(k-1)$-dimensional simplicial complex $\FF$ we have
\begin{equation}
\label{old_lower_bound}
\ex(n,\FF)\geq \ex_k(n,\FF^{k})+\sum_{r=0}^{k-1}\binom{n}{r},	
\end{equation}
where $\ex_k(n,H)$ denotes the Tur\'an number of a $k$-uniform hypergraph $H$. 
However, this lower bound could be arbitrarily far from the actual value, as could be seen for a $1$-dimensional simplicial complex $\FF$ with $\FF^2$ being a clique of order $t>3$. For this simplicial complex, it was shown in~\cite{CPS} that $\ex(n,\FF)= \Theta(n^{t-1})$, while $\ex_2(n,\FF^2)=\Theta(n^2)$.%\vspace{1mm}

It was further observed in~\cite{CPS} that, for some simplicial complex $\FF$, $\ex(n,\FF)$ is much closer to the so-called generalised Tur\'an numbers of related hypergraphs $H$.
More precisely, given $k$-uniform hypergraphs $T$ and $G$, let us denote by $\cN(T,G)$ the number of copies of $T$ in $G$. For a family $\cH$ of $k$-uniform hypergraphs, the \emph{generalised Tur\'an number} $\ex_k(n,T,\cH)$ is defined as
\[\ex_k(n,T,\cH)=\max\{\cN(T,G):\,G\text{ is an $n$-vertex $\cH$-free $k$-uniform hypergraph}\}.\]
In the case when $\cH=\{H\}$, we simply write $\ex_k(n,T,H)$ instead of $\ex_k(n,T,\{H\})$. The systematic study of generalised Tur\'an numbers was initiated by Alon and Shikhelman~\cite{ALS2016}, see also~\cite{GP2025} for a survey.
Conlon, Piga, and Sch\"ulke~\cite{CPS} showed that if $\FF$ is a simplicial complex whose maximal edges form a copy of $K_{t+1}^k$, that is, a complete $k$-uniform hypergraph on $t+1$ vertices, then
$\ex(n,\FF)=(1+o(1))\ex_k(n,K_t^k,K_{t+1}^k)$ for all $t\geq k\geq2$.%\vspace{1mm}

Given a $k$-uniform hypergraph $G$ and a subset $T\subseteq V(G)$, we say that $T$ forms a \emph{clique} in $G$ if $T$ is a singleton set, $T$ is a subset of some edge in $G$, or $|T|>k$ and $e\in E(G)$ for all $e\in\binom{T}{k}$.
For a family $\cH$ of $k$-uniform hypergraphs, we define $\ex_k^{\cl}(n,\cH)$ (resp. $\ex_k^{\cl+}(n,\cH)$) to be the maximum number of cliques (resp. cliques of order at least $k$) in an $n$-vertex $k$-uniform hypergraph $G$ that is $\cH$-free.
Generalising the observation~\eqref{old_lower_bound} in~\cite{CPS}, we show the following simple bounds on $\ex(n,\FF)$ in terms of the functions $\ex_k^{\cl}$ and $\ex_k^{\cl+}$.
For any hypergraph $H$ we denote by $\cD(H)$ the simplicial complex on the vertex set $V(H)$, whose edge set consists of all singletons and the downward closure of $E(H)$, i.e., all subsets of edges of $H$.

\begin{prop}
\label{LB}
Let $k\geq2$ be an integer and $\FF$ a $(k-1)$-dimensional simplicial complex. 
Let $\cH_{\FF}$ be the family of all $k$-uniform hypergraphs $H$ on the vertex set $V$, such that $\cD(H)$ contains a copy of $\FF$.
Then 
\begin{enumerate}[(1)]
\item $\ex_k^{\cl}(n,\cH_{\FF}) \leq \ex(n,\FF) \leq \ex_k^{\cl+}(n,\cH_{\FF})+\sum_{r=0}^{k-1}\binom{n}{r}$;
\item $\ex(n,\FF) \geq \max\left\{\ex_s^{\cl+}(n,\FF^{s})+\sum_{r=0}^{s-1}\binom{n}{r}:\,s=2,\dots,k\right\}$.
\end{enumerate}
\end{prop}

Note that~\Cref{LB} implies the aforementioned result from~\cite{CPS} on simplicial complexes whose maximal edges form a copy of $K^k_{t+1}$ with $t\geq k\geq2$. Indeed, for any such $\FF$ the unique minimal member in $\cH_{\FF}$ is a copy of $K^k_{t+1}$. Hence, by item (1) of~\Cref{LB} we have $\ex(n,\FF) \leq \ex_k^{\cl+}(n,K^{k}_{t+1})+O(n^{k-1})$. On the other hand, by item (2) of~\Cref{LB} we have $\ex(n,\FF) \geq \ex_k^{\cl+}(n,K^{k}_{t+1})+O(n^{k-1})$. Since there is no clique of order larger than $t$ in a $K^k_{t+1}$-free $k$-uniform hypergraph and the number of cliques of order smaller than $t$ is at most $O(n^{t-1})$, the asymptotics of $\ex(n,\FF)$ is given by $\ex_k(n,K_t^k,K_{t+1}^k)$. In the graph case (i.e., $k=2$), we can obtain an exact result. A theorem of Zykov~\cite{Z} shows that among all $n$-vertex $K_{t+1}$-free graphs, the Tur\'an graph $T(n,t)$ maximises the number of copies of $K_s$ for any $s\in[t]$. Therefore, $\ex_2(n,K_t,K_{t+1})=\cN(K_t,T(n,t))$.%\vspace{1mm}

Unfortunately, we do not have many other exact results implied by~\Cref{LB}. Generally speaking, if we have a $(k-1)$-dimensional simplicial complex $\FF$, such that $\FF=\cD(\FF^k)$ and $\ex_k^{\cl+}(n,\FF^k)$ is known, then one can apply~\Cref{LB} to obtain $\ex(n,\FF)$. To the best of our knowledge, the function $\ex_k^{\cl+}(n,\FF^k)$ is determined exactly only for $\FF^k$ being a graph star~\cite{CR} or path~\cite{CC}, and asymptotically for a larger class of trees~\cite{GMP}, book graphs, and fan graphs~\cite{G}.%\vspace{1mm}

By~\Cref{LB} we have that for any $(k-1)$-dimensional simplicial complex $\FF$, $\ex(n,\FF)\geq\ex_k^{\cl+}(n,\FF^{k})+\sum_{r=0}^{k-1}\binom{n}{r}$, which is also called the \emph{trivial lower bound} for $\ex(n,\FF)$. Simplicial complexes whose extremal number attains the trivial lower bound are of particular interest.
\begin{definition}
Let $k\geq2$ be an integer and $\FF$ a $(k-1)$-dimensional simplicial complex. We say that $\FF$ is trivial if for sufficiently large $n\in\NN$
\[\ex(n,\FF)=\ex_k^{\cl+}(n,\FF^{k})+\sum_{r=0}^{k-1}\binom{n}{r}.\]
Furthermore, the simplicial complex $\FF$ is asymptotically trivial if 
\[\ex(n,\FF)=\ex_k^{\cl+}(n,\FF^{k})+\sum_{r=0}^{k-1}\binom{n}{r}+o(n^{k-1}).\]
\end{definition}

The above definition naturally leads to the question of which simplicial complexes are trivial. Our first main result makes progress in this direction. We say that a hypergraph is \emph{edge-degenerate}, if there exists an ordering $e_1,\dots,e_m$ of its edges, such that $\forall\,i\in\{2,\dots,m\},\,\exists\,j\in\{1,\dots,i-1\}:\,e_{i}\cap\left(\bigcup_{r=1}^{i-1}e_r\right)\subseteq e_j$. A \emph{linear cycle} $C^k_t$ is a $k$-uniform hypergraph consisting of $t\geq3$ edges cyclically ordered, such that any two consecutive edges share exactly one vertex, any two non-consecutive edges are disjoint, and every vertex is contained in at most two edges.

\begin{thm}
\label{new:trivial}
Let $k\geq2$ be an integer. Let $\FF$ be a $(k-1)$-dimensional simplicial complex and $\cE(\FF)$ the set of maximal edges of $\FF$. Then $\FF$ is trivial, if one of the following holds:
\begin{enumerate}[(i)]
\item\label{(i)} $\cE(\FF)=E(\FF^k)$;
\item\label{(ii)} $\cE(\FF)$ consists of pairwise disjoint edges;
\item\label{(iii)} $\lvert\cE(\FF)\rvert\leq2$;
\item\label{(iv)} $k\geq3$ and $\cE(\FF)=\left\{\{1,2,\dots,k\},\{1,k+1,\dots,2k-2\},\{2,k+1,\dots,2k-2\}\right\}$.
\end{enumerate}
Moreover,
\begin{enumerate}
\item if $\FF$ satisfies~\eqref{(i)} and $\FF^k$ is edge-degenerate or, when $k\geq3$, a linear cycle, then $\ex(n, \FF)=\Theta(n^{k-1})$;
\item if $\FF$ satisfies~\eqref{(ii)}, then for sufficiently large $n\in\NN$ we have
\[\ex(n,\FF)=\sum_{r=1}^{t-1}\sum_{i=1}^{r}\binom{t-1}{r}\binom{n-t+1}{k-i}+\sum_{r=0}^{k-1}\binom{n}{r},\]
where $t$ denotes the number of edges in $\FF^k$.
\end{enumerate}
\end{thm}

\Cref{new:trivial}~\eqref{(i)} yields that $\ex(n,\cD(C^k_t))=\Theta(n^{k-1})$ for every $k\geq3$. The condition $k\geq3$ is necessary, as illustrated by the simplicial complex $\cD(C^2_4)$. Indeed, since $\ex_2(n,C^2_4)=\Theta(n^{3/2})$~\cite{B66,R58} and a straightforward argument gives $\ex_2^{\cl+}(n,C^2_4)=\Theta(\ex_2(n,C^2_4))$, we have that $\ex(n,\cD(C^2_4))=\Theta(n^{3/2})$. Let us also remark that~\Cref{new:trivial}~\eqref{(iii)} is best possible, as we shall see in~\Cref{thm:jump} that for every $k\geq2$ there are $(k-1)$-dimensional simplicial complexes $\FF$ with $\lvert\cE(\FF)\rvert=3$, such that $\ex(n,\FF)$ is much larger than the trivial lower bound.

In the next result, we provide a large class of simplicial complexes that are asymptotically trivial. Let $M^k_t$ denote a $k$-uniform matching on $t$ edges, and recall that $C^k_{2t}$ denotes a $k$-uniform linear cycle on $2t$ edges.

\begin{thm}
\label{new:asymptrivial}
Let $k\geq3$ and $t\geq2$ be integers. Let $\FF$ be a $(k-1)$-dimensional simplicial complex, such that $\FF$ contains a copy of $M^k_t$ and is contained in a copy of $\cD(C^k_{2t})$. Then $\FF$ is asymptotically trivial and we have
 \[\ex(n,\FF)=2^{t-1}\binom{n}{k-1}+o(n^{k-1}).\]
\end{thm}

Note that~\Cref{new:asymptrivial} extends a result of Conlon, Piga, and Sch\"ulke~\cite[Theorem~1.2]{CPS}, who considered the special case in which $\cE(\FF)=\left\{\{1,2,3\},\{4,5,6\},\{1,6\}\right\}$. For this special case, they obtained a sharper estimate of $\ex(n,\FF)$ up to an additive constant.

In contrast to the trivial and asymptotically trivial cases considered above, the next result provides several non-trivial $2$-dimensional simplicial complexes and bounds on their extremal numbers.
\begin{prop}
\label{thm:non-trivial}
Let $\FF_1,\FF_2,\FF_3,\FF_4$ be simplicial complexes described by their maximal edges as below. Then all of them are non-trivial and their extremal numbers satisfy the following:\\
\scalebox{0.93}[1]{
\begin{tabular}{ll}
$\cE(\FF_1)=\left\{\{1,2,3\},\{1,4,5\},\{3,4\},\{3,5\}\right\}$, & $\ex(n,\FF_1)=n^3/27+o(n^3)$;\\
$\cE(\FF_2)=\left\{\{1,2,3\},\{4,5,6\},\{1,4\},\{1,5\},\{1,6\}\right\}$, & $\ex(n,\FF_2)=n^3/27+o(n^3)$;\\
$\cE(\FF_3)=\left\{\{1,2,3\},\{4,5,6\},\{1,4\},\{1,5\},\{2,4\},\{2,5\}\right\}$, & $\ex(n,\FF_3)=n^3/27+o(n^3)$;\\
$\cE(\FF_4)=\left\{\{1,2,3\},\{4,5,6\},\{1,4\},\{1,6\},\{2,4\},\{2,5\},\{3,5\},\{3,6\}\right\}$, & $\Omega(n^{5/2})=\ex(n,\FF_4)=O(n^{25/8})$.
\end{tabular}
}
\end{prop}

We conclude by examining another aspect of the extremal numbers of simplicial complexes. In Tur\'an-type problems, a natural question is how the extremal number changes when an edge is added to the forbidden graph or hypergraph. Here, we address the analogous question for simplicial complexes and obtain the following result.

\begin{thm}
\label{thm:jump}
Let $k$ and $t$ be integers with $k\geq t\geq 2$. Let $\FF$ and $\HH$ be $(k-1)$-dimensional simplicial complexes, where $\cE(\FF)$ is obtained from $\cE(\HH)$ by adding a $t$-edge. Then
\[\frac{\ex(n,\FF)}{\ex(n,\HH)}=O(n^{t-1}).\]
Moreover, for every $k\geq t\geq2$ there are such $\FF$ and $\HH$ with $|\cE(\HH)|=t$ and $\frac{\ex(n,\FF)}{\ex(n,\HH)}=\Theta(n^{t-1})$.
\end{thm}

Note that in the case $t=2$,~\Cref{thm:jump} implies that there exists a $(k-1)$-dimensional simplicial complex $\FF$, such that $\lvert\cE(\FF)\rvert=3$ and $\ex(n,\FF)$ is much larger than its trivial lower bound. This shows in particular that~\Cref{new:trivial}~\eqref{(iii)} could not be extended to simplicial complexes with more than two maximal edges.

%\vspace{1mm}
The paper is structured as follows. We first collect some notation and definitions used throughout the paper in~\Cref{notations_and_definitions}. In~\Cref{preliminary} we prove~\Cref{LB} and several auxiliary results. In particular, we prove a generalised hypergraph Tur\'an result (\Cref{cliquelingen}) in~\Cref{preliminary} that is interesting on its own.
Then, we proceed to prove~\Cref{new:trivial} in~\Cref{trivial},~\Cref{new:asymptrivial} in~\Cref{asymptotically-trivial},~\Cref{thm:non-trivial} in~\Cref{non-trivial}, and~\Cref{thm:jump} in~\Cref{adding_edge}. Finally, we state some concluding remarks in~\Cref{conclusions}, including a connection of simplicial Tur\'an problems to the so-called Berge hypergraphs and a variety of open questions.

%%%%%%%%%%%%%%%%%%%%%%%%%%%%%%%%%%%%%%%%%%%%%%%%%%%%%%%%%%%%%%%%%%%%%%
%%%%%%%%%%%%%%%%%%%%%%%%%%%%%%%%%%%%%%%%%%%%%%%%%%%%%%%%%%%%%%%%%%%%%%
\section{Notations and definitions}
\label{notations_and_definitions}
%%%%%%%%%%%%%%%%%%%%%%%%%%%%%%%%%%%%%%%%%%%%%%%%%%%%%%%%%%%%%%%%%%%%%%
%%%%%%%%%%%%%%%%%%%%%%%%%%%%%%%%%%%%%%%%%%%%%%%%%%%%%%%%%%%%%%%%%%%%%%
Let $k\geq2$ be an integer. Given a $(k-1)$-dimensional simplicial complex $\FF$ and $r\in\{0,1,\dots,k\}$, we let $\FF^{r}$ denote the $r$-uniform hypergraph on the vertex set $V(\FF)$ consisting of all $r$-element edges (abbreviated as \emph{$r$-edges}) in $\FF$. We also call $\FF^r$ the \emph{$r$-th layer} of $\FF$.%\vspace{1mm}

The \emph{generating set} of $\FF$, denoted by $\cE(\FF)$, is the set of maximal edges in $\FF$. Since every $\FF$ is uniquely determined by $\cE(\FF)$, we often identify $\FF$ with its generating set $\cE(\FF)$.%\vspace{1mm}

For any hypergraph $H$, recall that $\cD(H)$ denotes the simplicial complex on the vertex set $V(H)$, whose edge set consists of all singletons and the downward closure of $E(H)$.%\vspace{1mm}

Given $k$-uniform hypergraphs $T$ and $G$, $\cN(T,G)$ denotes the number of copies of $T$ contained in $G$. Moreover, let $\cN^{\cl+}(G)$ denote the number of cliques of order at least $k$ contained in $G$.%\vspace{1mm}

Lastly, we introduce (or recall) the definitions of some specific hypergraphs that will be frequently used in this paper. Let $t\geq2$ be an integer. We let $K_{t}^k$ denote a complete $k$-uniform hypergraph on $t$ vertices, in particular, when $k=2$ we simply write $K_t=K^2_t$. A $k$-uniform hypergraph is a \emph{linear cycle} of length $t+1$, denoted by $C^k_{t+1}$, if it consists of $t+1$ edges cyclically ordered, such that any two consecutive edges share exactly one vertex, any two non-consecutive edges are disjoint, and every vertex is contained in at most two edges. A $k$-uniform \emph{linear path} $P^k_t$ of length $t$ is obtained from $C^k_{t+1}$ by deleting an arbitrary edge. In particular, we denote by $C_t=C^2_t$ and $P_t=P^2_t$ the graph cycle and graph path on $t$ edges, respectively. A $k$-uniform hypergraph is a \emph{matching} of size $t$, denoted by $M^k_t$, if it consists of $t$ pairwise disjoint edges.

%%%%%%%%%%%%%%%%%%%%%%%%%%%%%%%%%%%%%%%%%%%%%%%%%%%%%%%%%%%%%%%%%%%%%%
%%%%%%%%%%%%%%%%%%%%%%%%%%%%%%%%%%%%%%%%%%%%%%%%%%%%%%%%%%%%%%%%%%%%%%
\section{Proof of~\Cref{LB} and preliminary results}
\label{preliminary}
%%%%%%%%%%%%%%%%%%%%%%%%%%%%%%%%%%%%%%%%%%%%%%%%%%%%%%%%%%%%%%%%%%%%%%
%%%%%%%%%%%%%%%%%%%%%%%%%%%%%%%%%%%%%%%%%%%%%%%%%%%%%%%%%%%%%%%%%%%%%%

%%%%%%%%%%%%%%%%%%%%%%%%%%%%%%%%%%%%%%%%%%%%%%%%%%%%%%%%%%%%%%%%%%%%%%
\subsection{Proof of~\Cref{LB}}
%%%%%%%%%%%%%%%%%%%%%%%%%%%%%%%%%%%%%%%%%%%%%%%%%%%%%%%%%%%%%%%%%%%%%%
We first prove item (1). For the lower bound, we consider an $\cH_{\FF}$-free $k$-uniform hypergraph $G$ with $\ex_k^{\cl}(n,\cH_{\FF})$ cliques. Let $\HH$ be a simplicial complex whose edge set consists of all the cliques in $G$. Note that $\HH^k=G$. Suppose that $\HH$ contains a copy of $\FF$. Then for each $r<k$, every $r$-edge in that copy of $\FF$ is either a singleton set or a subset of some $k$-edge in $\HH$, otherwise it is not a clique in $G$. But this means that $\cD(G)=\cD(\HH^k)$ contains a copy of $\FF$, thus $G\in\cH_{\FF}$, a contradiction.%\vspace{1mm}

For the upper bound in item (1), let $\HH$ be an $n$-vertex $\FF$-free simplicial complex. Observe that $\HH^k$ is $\cH_{\FF}$-free by the definition of $\cH_{\FF}$. Moreover, each edge in $\HH$ of size larger than $k$ is a clique in $\HH^k$. This implies that the number of edges in $\HH$ of size at least $k$ is at most $\ex_k^{\cl+}(n,\cH_{\FF})$. The number of edges of size smaller than $k$ is at most $\sum_{r=0}^{k-1}\binom{n}{r}$, completing the proof of item (1).%\vspace{1mm}

We proceed to prove item (2). Fix any $s\in\{2,\dots,k\}$. Let $G$ be an $n$-vertex $\FF^s$-free $s$-uniform hypergraph with $\cN^{\cl+}(G)=\ex_s^{\cl+}(n,\FF^{s})$. Let $\HH$ be a hypergraph on the vertex set $V(G)$, whose edges are the cliques in $G$ and additionally all subsets of vertices of size at most $s-1$. Then $\HH$ contains $\ex_s^{\cl+}(n,\FF^{s})+\sum_{r=0}^{s-1}\binom{n}{r}$ edges. Observe that $\HH$ is a simplicial complex, as its edge set contains all singletons and is closed under taking subsets. Moreover, the $s$-edges of $\HH$ form an $\FF^s$-free hypergraph, thus $\HH$ contains no copy of $\FF$, completing the proof.\qed

%%%%%%%%%%%%%%%%%%%%%%%%%%%%%%%%%%%%%%%%%%%%%%%%%%%%%%%%%%%%%%%%%%%%%%
\subsection{Counting cliques in an $M^k_t$-free hypergraph}
%%%%%%%%%%%%%%%%%%%%%%%%%%%%%%%%%%%%%%%%%%%%%%%%%%%%%%%%%%%%%%%%%%%%%%
Given $k$-uniform hypergraphs $T$ and $H$, recall that $\ex_k(n,T,H)$ is the largest number of copies of $T$ in an $n$-vertex $H$-free $k$-uniform hypergraph. For sufficiently large $n$, Liu and Wang~\cite{LW} determined the value of $\ex_k(n,K^k_r,M_t^k)$ when $r\leq k+t-1$ and showed that $\ex_k(n,K^k_r,M_t^k)=O(n^{k-2})$ when $r\geq k+t$. Here we determine the value of $\ex_k^{\cl+}(n,M_t^k)$ for sufficiently large $n$.

\begin{prop}
\label{propmatch} 
Let $k,t\in\NN$ with $k\geq2$. For sufficiently large $n\in\NN$, we have
\begin{equation}
\label{matchclique}
\ex_k^{\cl+}(n,M_t^k)=\sum_{r=1}^{t-1}\sum_{i=1}^{r}\binom{t-1}{r}\binom{n-t+1}{k-i}=(2^{t-1}-1)\binom{n}{k-1}+o(n^{k-1}).
\end{equation}
\end{prop}

To prove~\Cref{propmatch}, we will use the following result by Gerbner~\cite{G2}.

\begin{lem}[{\cite[Proposition~2.1]{G2}}]
\label{coverin}
Let $k,t\geq2$ and $H$ be an $M_t^k$-free $k$-uniform hypergraph. Then there exists a set $A$ of at most $t-1$ vertices and a set $B$ of at most $k(2t-2)$ vertices such that every edge in $H$ contains at least one vertex from $A$ or at least two vertices from $B$. Moreover, if for every such pair $(A,B)$ we have $|A|=t-1$, then there is a pair where $B=\emptyset$.
\end{lem}

\begin{proof}[Proof of~\Cref{propmatch}]
When $t=1$, it holds trivially that both sides of~\eqref{matchclique} are equal to $0$.%\vspace{1mm} 

Assume that $t\geq2$. We start with the lower bound on $\ex_k^{\cl+}(n,M_t^k)$. Let $S^k_{n,t-1}$ be an $n$-vertex $k$-uniform hypergraph, whose edges are selected by fixing a set $A$ of $t-1$ vertices and taking every $k$-element set that intersects $A$. Observe that $S^k_{n,t-1}$ contains no copy of $M^k_t$, and each $T\subseteq V(S^k_{n,t-1})$ forms a clique if and only if $|T\setminus A|\leq k-1$. Then we have
$\ex_k^{\cl+}(n,M_t^k)\geq\cN^{\cl+}(S^k_{n,t-1})$ and 
\begin{equation}
\label{stars}
\cN^{\cl+}(S^k_{n,t-1})=\sum_{r=1}^{t-1}\sum_{i=1}^{r}\binom{t-1}{r}\binom{n-t+1}{k-i}.
\end{equation}

For the upper bound, let $H$ be an $n$-vertex $M_t^k$-free $k$-uniform hypergraph with $\cN^{\cl+}(H)=\ex_k^{\cl+}(n,M_t^k)$. We say that a pair $(A,B)\in2^{V(H)}\times2^{V(H)}$ is \emph{good} if $|A|\leq t-1$, $|B|\leq k(2t-2)$, and every edge in $H$ contains at least one vertex from $A$ or at least two vertices from $B$. By~\Cref{coverin} there exists at least one good pair $(A,B)$. If we can find a good pair $(A,B)$ with $|A|\leq t-2$, then we consider two types of cliques $T$ with $|T|\geq k$:
\begin{enumerate}[(a)]
    \item $|T\setminus A|\leq k-1$;
    \item $|T\setminus A|\geq k$.
\end{enumerate}
There are at most $\cN^{\cl+}(S^k_{n,t-2})$ cliques of type (a). On the other hand, if $T$ is of type (b), then $T$ contains at most $k-2$ vertices outside $A\cup B$. Indeed, otherwise we can take $k-1$ vertices from $T\setminus(A\cup B)$ and, since $|T\setminus A|\geq k$, another vertex from $T\setminus A$, which forms an edge in $H$ that is disjoint from $A$ and intersects $B$ in at most one vertex, a contradiction. This implies that there are $O(n^{k-2})$ cliques of type (b), namely, the total number of cliques of order at least $k$ is at most $\cN^{\cl+}(S^k_{n,t-2})+O(n^{k-2})<\cN^{\cl+}(S^k_{n,t-1})$. If for every good pair $(A,B)$ we always have $|A|=t-1$, then by~\Cref{coverin} we have one good pair $(A,B)$ with $B=\emptyset$. In this case $H$ is contained in a copy of $S^k_{n,t-1}$, completing the proof.
\end{proof}

%%%%%%%%%%%%%%%%%%%%%%%%%%%%%%%%%%%%%%%%%%%%%%%%%%%%%%%%%%%%%%%%%%%%%%
\subsection{Counting cliques in a hypergraph without linear path or cycle}
%%%%%%%%%%%%%%%%%%%%%%%%%%%%%%%%%%%%%%%%%%%%%%%%%%%%%%%%%%%%%%%%%%%%%%
In this subsection we prove the following generalised hypergraph Tur\'an result, which is of independent interest. Recall that for integers $k\geq2$ and $\ell\geq1$, we let $S^k_{n,\ell}$ be an $n$-vertex $k$-uniform hypergraph, consisting of all edges that intersect a fixed set of $\ell$ vertices.

\begin{thm}
\label{cliquelingen}
Let $k\geq3$ and $t\geq4$ be integers and $\ell=\lfloor(t-1)/2\rfloor$. Let $Q\in\{P^k_t,C^k_t\}$.
\begin{enumerate}[(i)]
	\item If $k\leq r\leq\ell+k-1$, then we have $\ex_k(n,K_r^k,Q)=(1+o(1))\cN(K^k_r,S_{n,\ell}^k)$.
	\item If $r\geq\ell+k$, then we have $\ex_k(n,K_r^k,Q)=o(n^{k-1})$.
\end{enumerate}
\end{thm}

We remark that when $r=k$, $\ex_k(n,K_r^k,Q)$ reduces to the ordinary Tur\'an number for linear paths or linear cycles. In this case, the exact values of $\ex_k(n,P^k_t)$ and $\ex_k(n,C^k_t)$ were determined by F\"uredi, Jiang, and Seiver~\cite{FJS14} and by F\"uredi and Jiang~\cite{FJ14}, respectively.

To prove~\Cref{cliquelingen} we will adopt an argument by Kostochka, Mubayi, and Verstra\"ete~\cite{KMV}. We say that a $k$-uniform hypergraph $G$ is \emph{$\ell$-full} if every $(k-1)$-element subset of $V(G)$ is contained in either none or at least $\ell$ edges. Extending a lemma from~\cite{KMV}, we prove the following.

\begin{lem}
\label{kmvlem1}
Let $k\geq2$ be an integer. If an $\ell$-vertex $k$-uniform hypergraph $Q$ is edge-degenerate or, when $k\geq3$, a linear cycle, then every non-empty $\ell$-full $k$-uniform hypergraph contains a copy of $Q$.
\end{lem}
\begin{proof}[Proof of~\Cref{kmvlem1}]
The case when $Q$ is a linear cycle was shown in~\cite[Lemma~3.2]{KMV}. Assume that $Q$ is an edge-degenerate $k$-uniform hypergraph on $\ell$ vertices. We shall prove the statement by induction on the number of edges in $Q$. When $Q$ has at most one edge, the statement follows trivially. When $Q$ has at least two edges, there are two distinct edges $e,e’\in E(Q)$ such that $Q-e$ is edge-degenerate and each vertex of $e$ either has degree $1$ in $Q$ or is contained in $e'$. Here $Q-e$ denotes the subhypergraph of $Q$ obtained by deleting the edge $e$. By our induction hypothesis, there is a copy of $Q-e$ in $H$. Let $A\subseteq V(H)$ be the set of vertices spanned by the edges of this copy. Let $S\subseteq A$ be the set corresponding to $e\cap e’$. Consider an edge $X$ of $H$ that contains $S$ and the largest number of vertices in $V(H)\setminus A$. Such an edge exists because in the worst case we can take $X$ to be the copy of $e'$. If $X\cap A=S$, then the copy of $Q-e$ and $X$ forms a copy of $Q$ with $X$ playing the role of $e$, as desired. If $X\cap A \neq S$, then we have $S\subseteq Y\subseteq X\setminus\{a\}$, for some $a\in X\cap A$ and $Y\subseteq X\setminus\{a\}$ with $|Y|=k-1$. Since $H$ is $\ell$-full for $\ell=|V(Q)|>|A|$, there is an edge $X’$ of $H$ containing $Y$ and one vertex in $V(H)\setminus A$, contradicting the maximality of $X$.
\end{proof}

For the proof of~\Cref{cliquelingen} we will also use the following result from~\cite{KMV}.
\begin{lem}[{\cite[Theorem~6.1]{KMV}}]
\label{kmvlem2}
Let $k\geq3$ and $t\geq4$ be integers and $\ell=\lfloor(t-1)/2\rfloor$. Let $G$ be an $n$-vertex $k$-uniform hypergraph that is $P^k_{t}$-free or $C^k_{t}$-free. If $G$ is $(\ell+1)$-full, then $|E(G)|=o(n^{k-1})$.
\end{lem}

\begin{proof}[Proof of~\Cref{cliquelingen}]
We first prove a general upper bound on $\ex_k(n,K^k_r,Q)$. Let $G$ be an $n$-vertex $Q$-free $k$-uniform hypergraph with $\cN(K^k_r,G)=\ex_k(n,K^k_r,Q)$. If $G$ is not $(\ell+1)$-full, then we apply the following procedure:
\begin{itemize}
    \item find a $(k-1)$-element subset $T\subseteq V(G)$ that is contained in at least one and at most $\ell$ edges;
    \item delete all the edges that contain $T$.
\end{itemize}
We iterate this procedure till the remaining hypergraph is $(\ell+1)$-full. This is feasible, as in the worst case we end up with an empty hypergraph and by definition it is $(\ell+1)$-full. This procedure takes at most $\binom{n}{k-1}$ iterations, since in each iteration we have picked a new $(k-1)$-element subset $T\subseteq V(G)$. Moreover, observe that when we delete all the edges containing $T$, we destroy all the cliques containing $T$. Since $T$ is contained in at most $\ell$ edges, there are at most $\ell$ vertices that can extend $T$ to a clique of order $r$. Namely, in each iteration we have destroyed at most $\binom{\ell}{r-k+1}$ cliques of order $r$. Therefore, by the end of this procedure we have destroyed at most
\[\binom{\ell}{r-k+1}\binom{n}{k-1}\] cliques of order $r$.%\vspace{1mm}

Let $H$ be the remaining hypergraph after the procedure. Since $H$ is $Q$-free and $(\ell+1)$-full, by~\Cref{kmvlem2} we have $|E(H)|=o(n^{k-1})$. If $H$ is empty, then there is no clique of order $r$ in $H$. Otherwise, since $Q$ is either edge-degenerate or a linear cycle, by~\Cref{kmvlem1} we know that $H$ is not $C$-full for some constant $C>0$ depending on $Q$. Now we can apply the above procedure by replacing $\ell+1$ with $C$. Note that this new procedure terminates only when all edges in $H$ are deleted, otherwise we would have a non-empty $C$-full hypergraph that is $Q$-free, a contradiction to~\Cref{kmvlem1}. Since $|E(H)|=o(n^{k-1})$, the number of $(k-1)$-element subsets $T\subseteq V(H)$ that are contained in some edges is at most $o(n^{k-1})$, namely, this new procedure takes at most $o(n^{k-1})$ iterations. Since in each iteration we destroy at most $2^{C}$ cliques of order $r$, the number of cliques in $H$ of order $r$ is at most $2^{C}o(n^{k-1})=o(n^{k-1})$. Therefore, for all $r\in\NN$ we have
\[\ex_k(n,K^k_r,Q)\leq\binom{\ell}{r-k+1}\binom{n}{k-1}+o(n^{k-1}).\]

If $k\leq r\leq\ell+k-1$, then
\[\cN(K^k_r,S^k_\ell)=\binom{\ell}{r-k+1}\binom{n}{k-1}+o(n^{k-1}),\]
which implies that $\ex_k(n,K^k_r,Q)\leq(1+o(1))\cN(K^k_r,S^k_\ell)$. On the other hand, since $S^k_\ell$ is $\{P^k_{t},C^k_t\}$-free, we have the lower bound
$\ex_k(n,K^k_r,Q)\geq\cN(K^k_r,S^k_\ell)$, completing the proof of item (i).%\vspace{1mm}

If $r\geq\ell+k$, then $\binom{\ell}{r-k+1}=0$. Thus, $\ex_k(n,K^k_r,Q)= o(n^{k-1})$, proving item (ii).
\end{proof}

\Cref{cliquelingen} and (\ref{stars}) immediately yield the following, which we will use to prove~\Cref{new:asymptrivial}.
\begin{cor}
\label{cliquelin}
For integers $k\geq3$ and $t\geq2$, we have \[\ex_k^{\cl+}(n,C^k_{2t})=(2^{t-1}-1)\binom{n}{k-1}+o(n^{k-1}).\]
\end{cor}

%%%%%%%%%%%%%%%%%%%%%%%%%%%%%%%%%%%%%%%%%%%%%%%%%%%%%%%%%%%%%%%%%%%%%%
%%%%%%%%%%%%%%%%%%%%%%%%%%%%%%%%%%%%%%%%%%%%%%%%%%%%%%%%%%%%%%%%%%%%%%
\section{Trivial simplicial complexes}
\label{trivial}
%%%%%%%%%%%%%%%%%%%%%%%%%%%%%%%%%%%%%%%%%%%%%%%%%%%%%%%%%%%%%%%%%%%%%%
%%%%%%%%%%%%%%%%%%%%%%%%%%%%%%%%%%%%%%%%%%%%%%%%%%%%%%%%%%%%%%%%%%%%%%
In this section we prove~\Cref{new:trivial}.

%%%%%%%%%%%%%%%%%%%%%%%%%%%%%%%%%%%%%%%%%%%%%%%%%%%%%%%%%%%%%%%%%%%%%%
\begin{proof}[Proof of~\Cref{new:trivial} case~\eqref{(i)}]
%%%%%%%%%%%%%%%%%%%%%%%%%%%%%%%%%%%%%%%%%%%%%%%%%%%%%%%%%%%%%%%%%%%%%%
Let $k\geq2$ and $\FF$ be a $(k-1)$-dimensional simplicial complex satisfying $\cE(\FF)=E(\FF^k)$, equivalently, $\FF=\cD(\FF^k)$.%\vspace{1mm}

Let $\HH$ be an $n$-vertex $\FF$-free simplicial complex with the maximum number of edges. Observe that the edges in $\HH$ of size at least $k$ corresponds to the cliques in $\HH^k$ of order at least $k$. Moreover, $\HH^k$ must be $\FF^k$-free, otherwise $\cD(\HH^k)\subseteq \HH$ contains a copy of $\cD(\FF^k)=\FF$, a contradiction. Thus $\HH$ contains at most $\ex_k^{\cl+}(n,\FF^{k})$ edges of size $k$ or greater. As the number of edges of size smaller than $k$ is at most $\sum_{r=0}^{k-1}\binom{n}{r}$, we obtain $|E(\HH)|\leq\ex_k^{\cl+}(n,\FF^{k})+\sum_{r=0}^{k-1}\binom{n}{r}$, namely, $\FF$ is trivial.%\vspace{1mm}

Further, assume that $\FF^k$ is edge-degenerate or, when $k\geq3$, a linear cycle. We shall prove that $\ex(n,\FF)=\Theta(n^{k-1})$. Since $\FF$ is trivial, it suffices to show that $\ex_k^{\cl+}(n,\FF^k)=O(n^{k-1})$. Take an $n$-vertex $k$-uniform hypergraph $H$ that is $\FF^k$-free and such that $\cN^{\cl+}(H)=\ex_k^{\cl+}(n,\FF^k)$. If the hypergraph $H$ is empty, then $\cN^{\cl+}(H)=0$. If $H$ is non-empty, then by~\Cref{kmvlem1} there exists a constant $C>0$ depending only on $\FF^k$, such that $H$ is not $C$-full. Now that $H$ is not $C$-full, we can iteratively choose a $(k-1)$-element subset $S\subseteq V(H)$ that is contained in at least one and at most $C-1$ edges, and then we delete all the edges containing $S$. This procedure takes at most $\binom{n}{k-1}$ steps, and in each step we destroy at most $2^{C-1}$ cliques of order at least $k$. Since every non-empty subhypergraph of $H$ is not $C$-full, this procedure stops when the remaining hypergraph is empty. Therefore, $\cN^{\cl+}(H)\leq2^{C-1}\binom{n}{k-1}=O(n^{k-1})$. This completes the proof of case~\eqref{(i)}.
\end{proof}

%\vspace{1.3em}
%%%%%%%%%%%%%%%%%%%%%%%%%%%%%%%%%%%%%%%%%%%%%%%%%%%%%%%%%%%%%%%%%%%%%%
\begin{proof}[Proof of~\Cref{new:trivial} case~\eqref{(ii)}]
%%%%%%%%%%%%%%%%%%%%%%%%%%%%%%%%%%%%%%%%%%%%%%%%%%%%%%%%%%%%%%%%%%%%%%
Let $\FF$ be a $(k-1)$-dimensional simplicial complex, where $\cE(\FF)$ consists of pairwise disjoint edges. Let $t\geq1$ be the number of $k$-edges in $\cE(\FF)$. Since $\FF^k$ is a copy of $M^k_t$, we want to show that \[\ex(n,\FF) \leq \ex_k^{\cl+}(n,M^k_t)+\sum_{r=0}^{k-1}\binom{n}{r}\]
for sufficiently large $n$.
Let $\HH$ be an $\FF$-free simplicial complex on $n$ vertices. We first consider the case when $\HH^k$ is $M_t^k$-free. Then, since the number of edges in $\HH$ of size at least $k$ is at most $\ex_k^{\cl+}(n,M^k_t)$, we have that $\lvert E(\HH)\rvert\leq\ex_k^{\cl+}(n,M^k_t)+\sum_{r=0}^{k-1}\binom{n}{r}$. Now consider the case when $\HH^k$ contains a copy of $M_t^k$. In this case we must have that $\cE(\FF)\neq E(\FF^k)$, otherwise $\HH$ would contain a copy of $\FF$, a contradiction. Let $p\geq1$ be the number of edges in $\cE(\FF)$ of size at most $k-1$. Let $A_1,\dots,A_t$ be the edges in that copy of $M_t^k$. Let $\{B_1,\dots,B_q\}$ be a largest set of pairwise disjoint edges in $\HH^{k-1}$, where each of them is disjoint from $A=A_1\cup\dots\cup A_t$. Since $\HH$ is $\FF$-free, it holds that $q<p$. Let $B=B_1\cup\dots\cup B_q$. Due to the maximality of $q$, every $(k-1)$-edge in $\HH$ intersects $A\cup B$. This implies that every edge in $\HH$ has at most $k-2$ vertices outside $A\cup B$. Since $|A\cup B|=O(1)$, we have that there are $O(n^{k-2})$ edges in $\HH$, which is smaller than the trivial lower bound on $\ex(n,\FF)$. Plugging in the value of $\ex_k^{\cl+}(n,M^k_t)$ from~\Cref{propmatch}, we complete the proof.
\end{proof}

%\vspace{1.3em}
%%%%%%%%%%%%%%%%%%%%%%%%%%%%%%%%%%%%%%%%%%%%%%%%%%%%%%%%%%%%%%%%%%%%%%
\begin{proof}[Proof of~\Cref{new:trivial} case~\eqref{(iii)}]
%%%%%%%%%%%%%%%%%%%%%%%%%%%%%%%%%%%%%%%%%%%%%%%%%%%%%%%%%%%%%%%%%%%%%%
Let $\lvert\cE(\FF)\rvert\leq2$. If $\cE(\FF)$ consists of only $k$-edges or disjoint edges, then by~\Cref{new:trivial} case~\eqref{(i)} or~\eqref{(ii)} we are done. Therefore, we assume that $\cE(\FF)$ consists of two intersecting edges, one of size $k$ and another of size at most $k-1$. Let $t$ be the intersection size of these two edges, $1\leq t\leq k-2$. Let $\HH$ be an $n$-vertex $\FF$-free simplicial complex with the maximum number of edges. 
Since $|V(\FF)|\leq 2k-2$, the largest size of an edge in $\HH$ is at most $|V(\FF)|-1\leq 2k-3$.
As $\FF^k$ consists of a single edge, $\ex_k^{\cl+}(n,\FF^k)=0$. Namely, it suffices to show that $|E(\HH)| \leq\sum_{r=0}^{k-1}\binom{n}{r}$ for sufficiently large $n$.%\vspace{1mm}

We shall first consider the edges in $\HH$ of size at least $k-1$. For every $t$-element subset $B$ of $V=V(\HH)$, let $x(B)$ be the number of edges in $\HH$ of size at least $k-1$ that contain $B$. Then the number of edges in $\HH$ of size at least $k-1$ is upper bounded by $\sum_{B\in\binom{V}{t}}x(B)/\binom{k-1}{t}$.%\vspace{0.3mm}

For each $B\in\binom{V}{t}$, if $B$ is contained in a $k$-edge $A\in E(\HH)$, then for any edge $A’\in E(\HH)$ of size at least $k-1$ containing $B$, we have that $|A’\setminus A|\leq k-t-2$, otherwise there is a copy of $\FF$. In this case it holds that $x(B)=O(n^{k-t-2})$. If $B$ is not contained in any $k$-edge of $\HH$, then $x(B)$ counts the number of $(k-1)$-edges containing $B$, namely, $x(B)\leq\binom{n-t}{k-1-t}$. 
In either case, for sufficiently large $n$, $x(B)\leq\binom{n-t}{k-t-1}$. This 
implies that the number of edges of size at least $k-1$ in $\HH$ is at most $\binom{n}{t}\binom{n-t}{k-t-1}/\binom{k-1}{t} = \binom{n}{k-1}$. Therefore, $|E(\HH)| \leq \sum_{r=0}^{k-2} \binom{n}{r} + \binom{n}{k-1}$.
\end{proof}

%\vspace{1.7em}
For the proof of~\Cref{new:trivial} case~\eqref{(iv)}, we shall use the following result by Frankl and F\"uredi~\cite{FF85}.
\begin{lem}[{\cite[Theorem~2.1]{FF85}}]
\label{intersecting}
For every integer $k\geq3$ there exists a constant $c_k>0$ such that the following holds. Let $t_1,t_2\in\NN$ with $t_1+t_2<k$ and $\cF$ be a family of $k$-element subsets of a ground set of $n$ elements. If $|A\cap B|\notin\{t_1,t_1+1,\dots,k-t_2-1\}$ for all $A,B\in\cF$, then $|\cF|\leq c_k n^{\max\{t_1,t_2\}}$.
\end{lem}
%%%%%%%%%%%%%%%%%%%%%%%%%%%%%%%%%%%%%%%%%%%%%%%%%%%%%
\begin{proof}[Proof of~\Cref{new:trivial} case~\eqref{(iv)}]
%%%%%%%%%%%%%%%%%%%%%%%%%%%%%%%%%%%%%%%%%%%%%%%%%%%%%
The proof idea of this case is similar to the previous one, with a more careful counting argument. Let $k\geq3$ and $\FF$ be a $(k-1)$-dimensional simplicial complex with \[\cE(\FF)=\left\{\{1,2,\dots,k\},\{1,k+1,\dots,2k-2\},\{2,k+1,\dots,2k-2\}\right\}.\]
Let $\HH$ be an $n$-vertex $\FF$-free simplicial complex with the maximum number of edges. As $\FF^k$ consists of a single edge, $\ex_k^{cl+}(n,\FF^k)=0$. Our goal is to show that $|E(\HH)| \leq\sum_{r=0}^{k-1}\binom{n}{r}$. We first prove the case when $k\geq4$, the proof for $k=3$ is slightly different. We call a set of $t$ vertices a \emph{$t$-non-edge}, or simply \emph{non-edge}, if it is not an edge in $\HH$.%\vspace{0.7em}

\noindent
\textbf{Case 1:} $k\geq4$.%\vspace{1mm}

Since $|V(\FF)|=2k-2$, the largest size of an edge in $\HH$ is at most $|V(\FF)|-1=2k-3$. %\vspace{1mm} 

Observe that every two $(k+s)$-edges in $\HH$, $0\leq s\leq k-3$, intersect in either at most one vertex or at least $s+3$ vertices.
Indeed, suppose there are two edges $\{1, \dots, i, x_{1}, \dots, x_{k+s-i}\}$ and $\{1, \dots, i, y_{1}, \dots, y_{k+s-i}\}$ in $\HH$, where $2\leq i\leq s+2$ and $x_p\neq y_q$ for all $p, q \in [k+s-i]$. Then $\{1, 2, x_1, \dots, x_{k-2}\}, \{1, y_1, \dots, y_{k-2}\}, \{2, y_1, \dots, y_{k-2}\}\in E(\HH)$ would give a copy of $\cE(\FF)$, a contradiction. Thus, by applying~\Cref{intersecting} on $E(\HH^{k+s})$ with $t_1=2$ and $t_2=k-3$, we obtain that $|E(\HH^{k+s})|= O(n^{k-2})$. 
Let $m$ be the number of edges in $\HH$ of size at least $k$. Then we have $m=\sum_{s=0}^{k-3}|E(\HH^{k+s})|= O(n^{k-2})$.%\vspace{1mm}

To show that $\FF$ is trivial, we shall prove that the number of $(k-1)$-non-edges in $\HH$ is at least $m$. Let $A$ be an edge of size $k+s$, $0\leq s\leq k-3$. For each set $T\subseteq V(\HH)\setminus A$ with $|T|=k-2$, there are at least $|A|-1\geq k-1$ vertices $v\in A$ such that $T\cup \{v\}$ is a non-edge. Thus the total number of $(k-1)$-non-edges in $\HH$ is at least $m\binom{n-(2k-3)}{k-2}(k-1) /x$, where $x$ is the overcount factor, i.e., the largest number of edges of size at least $k$ in $\HH$ that intersect a given $(k-1)$-element set in exactly one vertex. Consider a fixed $(k-1)$-element set $T\subseteq V(\HH)$ and a fixed vertex $v\in T$. For $0\leq s\leq k-3$, let $X(T,v,s)$ be the set of all $(k+s)$-edges $e\in E(\HH)$ such that $e\cap T=\{v\}$. Recall that any two $(k+s)$-edges in $\HH$ intersect in either at most $1$ vertex or at least $s+3$ vertices. Then $X(T,v,s)^{-}=\{e\setminus\{v\}:\,e\in X(T,v,s)\}$ is a family of $(k+s-1)$-element sets, whose pairwise intersection size is either $0$ or at least $s+2$. Thus we have $|X(T,v,s)|=|X(T,v,s)^{-}|= O(n^{k-3})$ by~\Cref{intersecting} applied to $X(T,v,s)^{-}$ with $t_1=1$ and $t_2=k-3$. Hence, the number of edges in $\HH$ of size at least $k$ that intersect $T$ in $v$ is $\sum_{s=0}^{k-3}|X(T,v,s)|= O(n^{k-3})$, i.e., $x= O(n^{k-3})$. 
Therefore, the number of $(k-1)$-non-edges in $\HH$ is at least $m\binom{n-(2k-3)}{k-2}(k-1) /x \geq m\Omega(n)$ and thus for sufficiently large $n$
\[|E(\HH)| \leq \sum_{r=0}^{k-2} \binom{n}{r} + \left(\binom{n}{k-1} -m\Omega(n)\right) + m \leq \sum_{r=0}^{k-1} \binom{n}{r}.\]

\noindent
\textbf{Case 2:} $k=3$.%\vspace{1mm}

In this case we have $\cE(\FF)=\left\{\{1,2,3\},\{1,4\},\{2,4\}\right\}$. Let $\HH$ be an $n$-vertex $\FF$-free simplicial complex with the maximum number of edges. Note that each edge in $\HH$ has size at most $|V(F)|-1=3$. Let $\ell$ be the number of non-edges in $\HH^2$. As $\ex(n,\FF)=\sum_{r=0}^3|E(\HH^r)|=\sum_{r=0}^2\binom{n}{2}-\ell+|E(\HH^3)|$, to prove that $\FF$ is trivial, it suffices to show that $\ell\geq|E(\HH^3)|$ for $n\geq4$. Let
\[\mathcal{J}=\left\{(e,f):\,e\in E(\HH^3)\text{ and $f$ is a non-edge in $\HH^2$ with $|e\cap f|=1$}\right\}.\]
We shall count the elements of $\mathcal{J}$ in two ways. First, for each $e\in E(\HH^3)$ and $v\in V(\HH)\setminus e$, there is at most one edge in $\HH^2$ containing $v$ and another vertex in $e$, as otherwise there would be a copy of $\FF$ in $\HH$. Consequently,
\[|\mathcal{J}|\geq|E(\HH^3)|\cdot(n-3)\cdot2.\]
On the other hand, every two edges in $\HH^3$ intersect in at most one vertex, otherwise there would be a copy of $\FF$ in $\HH$, a contradiction. Hence, for each non-edge $f$ in $\HH^2$ and $v\in f$, the number of $e\in E(\HH^3)$ with $e\cap f=\{v\}$ is at most $(n-2)/2$. Therefore,
\[|\mathcal{J}|\leq \ell\cdot2\cdot\frac{n-2}{2}=\ell\cdot(n-2).\]
Combining the lower and upper bounds on $|\mathcal{J}|$ yields $M\geq|E(\HH^3)|$ whenever $n\geq4$, as desired.
\end{proof}

%%%%%%%%%%%%%%%%%%%%%%%%%%%%%%%%%%%%%%%%%%%%%%%%%%%%%%%%%%%%%%%%%%%%%%
%%%%%%%%%%%%%%%%%%%%%%%%%%%%%%%%%%%%%%%%%%%%%%%%%%%%%%%%%%%%%%%%%%%%%%
\section{Asymptotically trivial simplicial complexes}
\label{asymptotically-trivial}
%%%%%%%%%%%%%%%%%%%%%%%%%%%%%%%%%%%%%%%%%%%%%%%%%%%%%%%%%%%%%%%%%%%%%%
%%%%%%%%%%%%%%%%%%%%%%%%%%%%%%%%%%%%%%%%%%%%%%%%%%%%%%%%%%%%%%%%%%%%%%
\begin{proof}[Proof of~\Cref{new:asymptrivial}]
Recall that $\cH_{\FF}$ denotes the family of all $k$-uniform hypergraphs $H$, such that $\cD(H)$ contains a copy of $\FF$. Since $C_{2t}^k\in\cH_{\FF}$ and $\FF^k$ contains a copy of $M^k_t$, by~\Cref{LB}, we have 
\[\ex(n,\FF)\geq \ex_k^{\cl+}(n,\FF^k) +\sum_{r=0}^{k-1}\binom{n}{r} \geq\ex_k^{\cl+}(n,M^{k}_t)+\sum_{r=0}^{k-1}\binom{n}{r}\]
and
\[\ex(n,\FF)\leq \ex_k^{\cl+}(n,\cH_{\FF})+ \sum_{r=0}^{k-1}\binom{n}{r} \leq \ex_k^{\cl+}(n,C^k_{2t})+\sum_{r=0}^{k-1}\binom{n}{r}.\]
By~\Cref{propmatch} and~\Cref{cliquelin} it holds that $\ex_k^{\cl+}(n,C^{k}_{2t})=\ex_k^{\cl+}(n,M^k_t)+o(n^{k-1})$, which further yields that $\ex(n,\FF)\leq \ex_k^{\cl+}(n,\FF^k) +\sum_{r=0}^{k-1}\binom{n}{r} + o(n^{k-1})$.
\end{proof}

%%%%%%%%%%%%%%%%%%%%%%%%%%%%%%%%%%%%%%%%%%%%%%%%%%%%%%%%%%%%%%%%%%%%%%%%%
%%%%%%%%%%%%%%%%%%%%%%%%%%%%%%%%%%%%%%%%%%%%%%%%%%%%%%%%%%%%%%%%%%%%%%%%%
\section{Non-trivial $2$-dimensional simplicial complexes}
\label{non-trivial}
%%%%%%%%%%%%%%%%%%%%%%%%%%%%%%%%%%%%%%%%%%%%%%%%%%%%%%%%%%%%%%%%%%%%%%%%%
%%%%%%%%%%%%%%%%%%%%%%%%%%%%%%%%%%%%%%%%%%%%%%%%%%%%%%%%%%%%%%%%%%%%%%%%%
In this section we prove~\Cref{thm:non-trivial}. First we introduce the definition of a hypergraph blow-up. Let $H$ be a $k$-uniform hypergraph and $t\in\NN$. The \emph{$t$-blow-up} of $H$, denoted $H(t)$, is a $k$-uniform hypergraph obtained from $H$ by replacing every vertex $v\in V(H)$ with an independent set $I_v$ of size $t$, such that $I_v$'s are pairwise disjoint, and replacing every edge $e=\{v_1,\dots,v_k\}\in E(H)$ by a complete $k$-partite $k$-uniform hypergraph with parts $I_{v_1},\dots,I_{v_k}$. It is well-known that the blow-up operation does not change the Tur\'an density (see, e.g.,~\cite{PK11}), namely, we have $\ex(n,H)=\ex(n,H(t))+o(n^k)$.%\vspace{1mm}

For the proof of~\Cref{thm:non-trivial}, we will use the following two results. The first one is a classical theorem by Zykov~\cite{Z}, which determines the maximum number of cliques in an $n$-vertex graph that is $K_{t+1}$-free. The second, due to Baber and Talbot~\cite{BT12}, asymptotically determines the extremal number of a specific $3$-uniform hypergraph.
\begin{lem}[{\cite[Theorem~15]{Z}}]
\label{Zykov}
Given $t\geq2$ and sufficiently large $n\in\NN$, among all $K_{t+1}$-free $n$-vertex graphs, the balanced complete $t$-partite graph maximises the number of cliques of any order. In particular, $\ex_2^{\cl+}(n,K_{t+1})=(n/t)^t+\Theta(n^{t-1})$.
\end{lem}

\begin{lem}[{\cite[Theorem~12]{BT12}}]
\label{density}
Let $H$ be a hypergraph with $V(H)=[6]$ and \[E(H)=\{\{1,2,3\},\{1,2,4\},\{3,4,5\},\{1,5,6\}\}.\] Then for any fixed $t\in\NN$, $\ex_{3}(n,H(t))=n^3/27+o(n^3)$.
\end{lem}

\begin{proof}[Proof of~\Cref{thm:non-trivial}]
We first bound $\ex(n, \FF_i)$ for $i=1,2,3,4$, and then show that $\FF_i$'s are all non-trivial. Given a simplicial complex $\HH$, we denote by $m_r(\HH)$ the number of $r$-edges in $\HH$ and we use the notation $m_{\geq r}(\HH)$ for the number of edges in $\HH$ of size at least $r$.%\vspace{1mm}

We start by showing that $\ex(n,\FF_1)=n^3/27+o(n^3)$. Since $\FF_1^2$ contains a copy of $K_4$, by~\Cref{LB} and~\Cref{Zykov} we obtain the lower bound $\ex(n,\FF_1)\geq\ex_2^{\cl}(n,\FF_1^2)\geq\ex_2^{\cl}(n,K_4)=n^3/27+o(n^3)$.
For the upper bound, we let $\HH$ be an $n$-vertex $\FF_1$-free simplicial complex with $|E(\HH)|=\ex(n,\FF_1)$. Since the number of edges in $\HH$ of size smaller than $3$ is at most $O(n^2)$, which is negligible, our goal is only to bound the number of edges of larger size. Let $m_r=m_r(\HH)$ and $m_{\geq r}= m_{\geq r}(\HH)$, $r\geq 1$. Since $\HH$ is $\FF_1$-free and $|V(\FF_1)|=5$, we have $m_{\geq5}=0$.
The intersection of any distinct $4$-edges $A$ and $B$ in $E(\HH)$ has size either $0$ or $1$, otherwise the downward closure of $\{A,B\}$ would contain a copy of $\cE(\FF_1)$, a contradiction. Then every pair of vertices is contained in at most one $4$-edge, implying that $\HH^4$ contains at most $\binom{n}{2}$ edges, i.e., $m_4= O(n^2)$.
To bound the number of $3$-edges, we look at the $3$-uniform hypergraph $\HH^3$. Let $H$ be the $3$-uniform hypergraph defined in~\Cref{density}. Observe that $\cD(H)$ contains a copy of $\FF_1$. Since $\HH$ is $\FF_1$-free, $\HH^3$ must be $H$-free, thus, $m_3=|E(\HH^3)|\leq\ex_{3}(n,H)=n^3/27+o(n^3)$.%\vspace{1.3em}

Next, we determine $\ex(n, \FF_2)$ asymptotically. As $\FF_2^2$ contains a copy of $K_4$, by the same argument as in the previous proof we have $\ex(n,\FF_2)\geq n^3/27+o(n^3)$.
Now let $\HH$ be an $\FF_2$-free $n$-vertex simplicial complex. We want to show that $|E(\HH)|\leq n^3/27+o(n^3)$. As the number of edges in $\HH$ of size smaller than $3$ is at most $O(n^2)$, it suffices to bound the number of edges of larger size. 
Let $m_r=m_r(\HH)$ and $m_{\geq r}= m_{\geq r}(\HH)$, $r\geq 1$.
Since $\HH$ is $\FF_2$-free and $|V(\FF_2)|=6$, we have $m_{\geq6}=0$. If there are two $5$-edges $A$ and $B$ with non-empty intersection, then the downward closure of $\{A,B\}$ would contain a copy of $\cE(\FF)$, a contradiction. Hence, the $5$-edges in $\HH$ must be pairwise vertex-disjoint, namely, $m_5\leq n/5$. Moreover, the intersection of any two $4$-edges $C$ and $D$ can have size either $0$ or $3$, otherwise we get a contradiction as the downward closure of $\{C,D\}$ contains a copy of $\cE(\FF_2)$. If we take a maximal collection $\mathcal{C}$ of pairwise disjoint $4$-edges, then every remaining $4$-edge must intersect exactly one member of $\mathcal{C}$ in exactly $3$ vertices. Moreover, two $4$-edges intersecting different members of $\mathcal{C}$ must have distinct fourth vertices. Thus, it follows that $m_4=O(n)$. Lastly, observe that $\cD(H(2))$ contains a copy of $\FF_2$, where $H$ is the hypergraph defined in~\Cref{density} and $H(2)$ denotes the $2$-blow-up of $H$. This implies that $\HH^3$ is $H(2)$-free. By~\Cref{density} we have $m_3\leq n^3/27+o(n^3)$.%\vspace{1.3em}

Now we show that $\ex(n,\FF_3)=n^3/27+o(n^3)$. The lower bound follows by the same argument as above. Let $\HH$ be an $n$-vertex $\FF_3$-free simplicial complex. It remains to show that $|E(\HH)|\leq n^3/27+o(n^3)$. Since the number of edges of size smaller than $3$ is at most $O(n^2)$ and $\HH$ contains no edge of size at least $|V(\FF_3)|=6$, it suffices to bound the number of $3$-edges, $4$-edges, and $5$-edges, denoted by $m_3$, $m_4$, and $m_5$, respectively.
Let $H$ be the hypergraph defined in~\Cref{density}. Observe that $\cD(H(2))$ contains a copy of $\FF_3$, which in turn yields that $m_3\leq\ex_3(n,H(2))=n^3/27+o(n^3)$.
Observe that the intersection size of any two $5$-edges is either $0$ or $1$, otherwise their downward closure contains a copy of $\cE(\FF_3)$, a contradiction. Then, as every pair of vertices can be contained in at most one $5$-edge, we have $m_5\leq\binom{n}{2}$.
To bound the number of $4$-edges, observe that every edge in $\HH^4$ contains a pair of vertices that is contained in at most $3$ edges in $\HH^3$. Indeed, suppose there is an edge $\{a,b,c,d\}\in E(\HH^4)$, where $\{a,b\}$ and $\{c,d\}$ are both contained in at least $4$ edges in $\HH^3$. Then there exist two distinct vertices $x,y\in V(\HH)\setminus\{a,b,c,d\}$, such that $\{a,b,x\},\{c,d,y\}\in E(\HH^3)$. Then $\cD(\{a,b,c,d\},\{a,b,x\},\{c,d,y\})$ gives a copy of $\FF_3$ in $\HH$, a contradiction. Now, since every edge in $\HH^4$ contains a pair of vertices that is contained in at most $3$ edges in $\HH^3$, we can bound the number of $4$-edges by counting the number of such pairs and the number of $4$-edges containing each pair. This gives that $m_4\leq\binom{n}{2}\binom{3}{2}=O(n^2)$. %\vspace{1.3em}

Finally, we prove that $\Omega(n^{5/2})=\ex(n,\FF_4)= O(n^{25/8})$. For the lower bound, we shall construct an $\FF_4$-free $n$-vertex simplicial complex $\HH$ such that $|E(\HH)|=\Omega(n^{5/2})$.
Let $V$ and $W$ be two disjoint sets of vertices with $|V|+|W|=n$. Recall that $C_t$ denotes a graph cycle of length $t$. Let $G$ be a $C_4$-free graph of maximum size on the vertex set $V$. Then $G$ must be connected and it is known~\cite{B66} that $|E(G)|=\Omega(|V|^{3/2})$. Let $\HH$ be a simplicial complex on the vertex set $V\dot\cup W$, where
\[\cE(\HH)=\left\{\{u,v,w\}:\,\{u,v\}\in E(G),w\in W\right\}.\]
By letting $|V|$ and $|W|$ be linear in $n$, we obtain $|E(\HH)|\ge |E(G)|\cdot|W|=\Omega(n^{5/2})$.
It remains to verify that $\HH$ is $\FF_4$-free. Suppose for a contradiction that $\HH$ contains a copy $\FF$ of $\FF_4$, with two disjoint $3$-edges $\{u,v,w\}$ and $\{u',v',w'\}$,
where $w, w'\in W$. Moreover, assume without loss of generality that $\{u,u'\}, \{v, v'\} \in E(\FF)$.
Since $\{u, v\}, \{u', v'\} \in E(\FF)$ are subsets of the edges $\{u,v,w\}$ and $\{u',v',w'\}$, we see that $\HH^2$ restricted to $V$ contains a copy of $C_4$, a contradiction. This shows the lower bound $\ex(n,\FF_4)=\Omega(n^{5/2})$.%\vspace{1mm}

For the upper bound, we let $\HH$ be an $n$-vertex $\FF_4$-free simplicial complex with the maximum number of edges.
Let $m_r=m_r(\HH)$. Since $m_1+m_2=O(n^2)$ and there are no edges of size $6$ or greater, it suffices to show that $m_3+m_4+m_5=O(n^{25/8})$. Let $K^{3}_{2,2,2}$ denote the complete $3$-partite $3$-uniform hypergraph with each vertex class of size $2$. Suppose $\HH^3$ contains a copy of $K^{3}_{2,2,2}$ on the vertex set $\{a_1,a_2\}\dot\cup\{b_1,b_2\}\dot\cup\{c_1,c_2\}$. Then
\[\left\{\{a_1,b_1,c_1\},\{c_2,a_2,b_2\},\{a_1,c_2\},\{a_1,b_2\},\{b_1,c_2\},\{b_1,a_2\},\{c_1,a_2\},\{c_1,b_2\}\right\}\subseteq E(\HH),\]
which implies that $\HH$ contains a copy of $\FF_4$, a contradiction. Since $\HH^3$ is $K^{3}_{2,2,2}$-free, by a classical result of Erd\H{o}s~\cite{Er64}, we have that $m_3= O(n^{11/4})$. To bound $m_4$, for any two vertices $x,y\in V(\HH)$ we let $G_{x,y}$ denote the graph with the vertex set $V(\HH)$ and the edge set
\[E(G_{x,y})=\left\{\{z,z'\}:\,\{x,z,z'\},\{y,z,z'\}\in E(\HH^3)\right\}.\]
In addition, let $G'_{x,y}$ denote the induced subgraph of $G_{x,y}$ on the vertex set
\[V(G'_{x,y})=\left\{z\in V(\HH):\,\{x,y,z\}\in E(\HH^3)\right\}.\]
Note that $m_4\leq\sum_{\{x,y\}\subseteq V(\HH)}|E(G'_{x,y})|$. Moreover, observe that for any distinct $x,y\in V(\HH)$, $G_{x,y}$ is $C_4$-free. Indeed, suppose there exists a copy of $C_4$ in $G_{x,y}$, with the edges $\{a,b\},\{b,c\},\{c,d\},\{d,a\}$. Then
\[\left\{\{x,a,b\},\{y,c,d\},\{x,c\},\{x,d\},\{a,y\},\{a,d\},\{b,y\},\{b,c\}\right\}\subseteq E(\HH),\]
implying that $\HH$ contains a copy of $\FF_4$, a contradiction. As $G_{x,y}$ (and hence $G'_{x,y}$) is $C_4$-free, it holds that $|E(G_{x,y})|\leq n^{3/2}$ and $|E(G'_{x,y})|\leq |V(G'_{x,y})|^{3/2}$. Since $m_3=O(n^{11/4})$ and $|V(G'_{x,y})|$ counts the number of $3$-edges containing $\{x,y\}$, we have that $\sum_{\{x,y\}\subseteq V(\HH)}|V(G'_{x,y})|\leq3m_3=O(n^{11/4})$. Then
\begin{align*}
m_4\leq\sum_{\{x,y\}\subseteq V(\HH)}|V(G'_{x,y})|^{3/2}&\leq\left(\sum_{\{x,y\}\subseteq V(\HH)}|V(G'_{x,y})|\right)^{1/2}\left(\sum_{\{x,y\}\subseteq V(\HH)}|V(G'_{x,y})|^2\right)^{1/2}\\
&\leq O(n^{11/8})\cdot\left(\sum_{\{x,y\}\subseteq V(\HH)}\binom{|V(G'_{x,y})|}{2}\right)^{1/2},
\end{align*}
where in the second last step we used Cauchy--Schwarz inequality. Further, by double-counting the number of copies of $\left\{\{x,y,z\},\{x,y,z'\}\right\}$ in $\HH^3$ we have
\[\sum_{\{x,y\}\subseteq V(\HH)}\binom{|V(G'_{x,y})|}{2}=\sum_{\{z,z'\}\subseteq V(\HH)}|E(G_{z,z'})|=O(n^{7/2}),\]
from which it follows that $m_4=O(n^{25/8})$. To bound the number of $5$-edges, observe that if there are two $5$-edges in $\HH$ with an intersection of size $4$, then the downward closure of these two $5$-edges contains a copy of $\FF_4$, a contradiction. Therefore, every $4$-edge is contained in at most one $5$-edge, yielding that $m_5\leq m_4=O(n^{25/8})$. This shows that $m_3+m_4+m_5=O(n^{25/8})$ and thus $\ex(n,\FF_4)=O(n^{25/8})$.%\vspace{1.3em}

Now we show that $\FF_1$, $\FF_2$, $\FF_3$, and $\FF_4$ are non-trivial. Recall that a $2$-dimensional simplicial complex $\FF$ is trivial if for sufficiently large $n$
\[\ex(n,\FF)=\ex_3^{\cl+}(n,\FF^3)+\sum_{r=0}^{2}\binom{n}{r}.\]
Since the $3$-uniform linear cycle $C^3_4$ contains a copy of $\FF_1^3$, by~\Cref{cliquelin} we have \[\ex_3^{\cl+}(n,\FF_1^3)\leq\ex_3^{\cl+}(n,C^3_4)=\Theta(n^2).\] Therefore, if $\FF_1$ were to be trivial, then we would have $\ex(n,\FF_1)=\Theta(n^2)$, which is not the case. Similarly, if the simplicial complexes $\FF_2$, $\FF_3$, and $\FF_4$ were trivial, then their extremal numbers would have order $\Theta(n^2)$. This is because their third layer is a copy of $M^3_2$ and by~\Cref{propmatch} we have $\ex_3^{\cl+}(n,M^3_2)=\Theta(n^2)$. But we have seen that $\ex(n,\FF_2),\ex(n,\FF_3),\ex(n,\FF_4)\gg n^2$.
\end{proof}

\section{Extremal numbers upon edge addition}
\label{adding_edge}
In this section we prove~\Cref{thm:jump}.
\begin{proof}[Proof of~\Cref{thm:jump}]
Let $\FF$ and $\HH$ be $(k-1)$-dimensional simplicial complexes such that $\cE(\FF)=\cE(\HH)\cup\{T\}$ with $|T|=t$. We shall prove that, for sufficiently large $n\in\mathbb N$,
\begin{equation}
\label{eq:adding_edge}
\ex(n,\FF)\leq\binom{n}{t-1}\ex(n,\HH)+\sum_{r=0}^{t-2}\binom{n}{r}.
\end{equation}

Let $\GG$ be an $n$-vertex $\FF$-free simplicial complex. For every $S\in E(\GG^{t-1})$, let $\GG_S$ be the simplicial complex on the edge set
\[E(\GG_S)=\left\{A\subseteq V(\GG)\setminus S:\,A\cup S\in E(\GG)\right\}.\]
We claim that $|E(\GG_S)|\leq\ex(n,\HH)$. Without loss of generality we assume that $|V(\GG_S)|>|V(\HH)|$, as otherwise $|E(\GG_S)|\leq2^{|V(\HH)|}\leq\ex(n,\HH)$ already holds for large $n$. Suppose $\GG_S$ is not $\HH$-free, and let $\phi:V(\HH)\to V(\GG_S)$ be an embedding of $\HH$ into $\GG_S$. If $T\cap V(\HH)\neq\emptyset$, we choose a vertex $x\in T\cap V(\HH)$ and map $T\setminus\{x\}$ bijectively to $S$. Then the subcomplex of $\GG$ induced by $S\cup\phi\left(V(\HH)\setminus\{x\}\right)$ would contain a copy of $\FF$, a contradiction. Now we consider the case when $T\cap V(\HH)=\emptyset$. Since $|V(\GG_S)|>|V(\HH)|$, $\GG_S$ contains a vertex $v$ outside $\phi\left(V(\HH)\right)$. Then by mapping $T$ bijectively to $S\cup\{v\}$, we obtain a copy of $\FF$ in $\GG$, a contradiction. Therefore, $\GG_S$ is $\HH$-free and thus $|E(\GG_S)|\leq\ex(n,\HH)$.

Since every edge of $\GG$ of size at least $t-1$ is counted in $E(\GG_S)$ for some $S\in E(\GG^{t-1})$, we have
\[|E(\GG)|\leq\sum_{S\in E(\GG^{t-1})}|E(\GG_S)|+\sum_{r=0}^{t-2}\binom{n}{r}\leq\binom{n}{t-1}\ex(n,\HH)+\sum_{r=0}^{t-2}\binom{n}{r},\]
which proves~\eqref{eq:adding_edge}.

To complete the proof of~\Cref{thm:jump}, it remains to construct such $\FF$ and $\HH$ with $|\cE(\HH)|=t$ and $\frac{\ex(n,\FF)}{\ex(n,\HH)}=\Omega(n^{t-1})$. Let $B=\{b_1,\dots,b_{k-1}\}$ and $C=\{c_1,\dots,c_t\}$ be two disjoint sets, and let
\[\cE(\HH)=\{B\cup\{c_i\}:\,i\in[t]\}\quad\text{and}\quad\cE(\FF)=\cE(\HH)\cup\{C\}.\]
Since $\cE(\HH)=E(\HH^k)$ and $\HH^k$ is edge-degenerate,~\Cref{new:trivial}~\eqref{(i)} gives that $\ex(n,\HH)=\Theta(n^{k-1})$. Moreover, as $\FF^2$ is a copy of $K_{k+t-1}$, by~\Cref{LB} and~\Cref{Zykov} we have
\[\ex(n,\FF)\geq\ex_2^{\cl}(n,\FF^2)=\Omega(n^{k+t-2}).\]
Therefore, $\frac{\ex(n,\FF)}{\ex(n,\HH)}=\Omega(n^{t-1})$ as required.
\end{proof}

%%%%%%%%%%%%%%%%%%%%%%%%%%%%%%%%%%%%%%%%%%%%%%%%%%%%%%%
%%%%%%%%%%%%%%%%%%%%%%%%%%%%%%%%%%%%%%%%%%%%%%%%%%%%%%%
\section{Concluding remarks}
\label{conclusions}
%%%%%%%%%%%%%%%%%%%%%%%%%%%%%%%%%%%%%%%%%%%%%%%%%%%%%%%
%%%%%%%%%%%%%%%%%%%%%%%%%%%%%%%%%%%%%%%%%%%%%%%%%%%%%%%
\subsection{Simplicial complexes attaining the trivial lower bound}
In this work we continue the systematic study of simplicial Tur\'an problems initiated by Conlon, Piga, and Sch\"ulke~\cite{CPS}. Let $k\geq2$ be an integer and $\FF$ a $(k-1)$-dimensional simplicial complex. We give a general lower bound on $\ex(n,\FF)$, stating \[\ex(n,\FF) \geq \ex_k^{\cl+}(n,\FF^{k})+\sum_{r=0}^{k-1}\binom{n}{r},\]
and investigate the simplicial complexes attaining this bound (up to an error term $o(n^{k-1})$) for sufficiently large $n$. Such simplicial complexes are called \emph{trivial} (\emph{asymptotically trivial}). We provide broad classes of trivial and asymptotically trivial simplicial complexes and determine the extremal number for some of them. It is however believable that there are many unrevealed trivial simplicial complexes. This motivates the following problem

\begin{prob}
\label{p1}
For any integer $k\geq3$, characterise all trivial $(k-1)$-dimensional simplicial complexes.
\end{prob}

Note that~\Cref{new:trivial} yields that any $1$-dimensional simplicial complex $\FF$ is trivial. Indeed, since the generating set $\cE(\FF)$ is the union of $E(\FF^2)$ and some singleton sets, for $n\geq|V(\FF)|$, an $n$-vertex simplicial complex $\HH$ contains a copy of $\FF$ if and only if $\HH$ contains a copy of $\cD(\FF^2)$. Then by~\Cref{new:trivial}~\eqref{(i)} we have $\ex(n,\FF)=\ex(n,\cD(\FF^2))=\ex_2^{\cl+}(n,\FF^2)+\sum_{r=0}^{1}\binom{n}{r}$.%\vspace{1mm}

A problem in a similar spirit was raised by Conlon, Piga, and Sch\"ulke~\cite{CPS}, and we restate it here.

\begin{prob}[{\cite[Problem~6.2]{CPS}}]
\label{p2}
Given an integer $k\geq3$, characterise the $(k-1)$-dimensional simplicial complexes $\FF$ with $\ex(n,\FF)=\ex_k(n,\FF^{k})+\sum_{r=0}^{k-1}\binom{n}{r}$.
\end{prob}

Although the ordinary Tur\'an number $\ex_k(n,\FF^k)$ is in most cases more accessible than the generalised variant $\ex_k^{\cl+}(n,\FF^{k})$, it is less closely related to the function $\ex(n,\FF)$. Thus, we expect~\Cref{p2} to be more difficult than~\Cref{p1}.

We record two observations relevant to~\Cref{p1} and~\Cref{p2}. The first one is immediate: If $\FF$ is contained in a $(k-1)$-dimensional simplicial complex $\HH$ and $E(\FF^k)=E(\HH^k)$, then we have $\ex^{\cl+}(n,\HH^k)+\sum_{r=0}^{k-1}\binom{n}{r}\leq\ex(n,\FF)\leq\ex(n,\HH)$. In particular, if $\HH$ is a simplicial complex sought by~\Cref{p1} or~\Cref{p2}, then so is $\FF$. The second observation is stated as follows.

\begin{obs}
\label{obs:disjoint}
Let $\FF$ and $\HH$ be simplicial complexes, and let $\FF\sqcup\HH$ denote their vertex-disjoint union. Then, for every $n\geq|V(\HH)|$,
\[\max\left\{\ex(n,\FF),\ex(n,\HH)\right\}\leq\ex(n,\FF\sqcup\HH)\leq\max\left\{2^{|V(\HH)|}\ex(n-|V(\HH)|,\FF),\ex(n,\HH)\right\}.\]
In particular, if $\ex(n,\FF)\ll\ex(n,\HH)$, then $\ex(n,\FF\sqcup\HH)=\ex(n,\HH)$ for sufficiently large $n$.
\end{obs}
\begin{proof}[Proof of~\Cref{obs:disjoint}]
The lower bound is immediate, since every $\FF$-free or $\HH$-free simplicial complex is $\FF\sqcup\HH$-free. For the upper bound, let $\bS$ be an $n$-vertex $\FF\sqcup\HH$-free simplicial complex. If $\bS$ is $\HH$-free, then $|E(\bS)|\leq\ex(n,\HH)$. Otherwise, we fix a copy of $\HH$ in $\bS$ on the vertex set $W$, and let $\bS'$ be the subcomplex of $\bS$ induced by $V(\bS)\setminus W$. Then $\bS'$ must be $\FF$-free, as otherwise there is a copy of $\FF\sqcup\HH$ in $\bS$, a contradiction. Since every $e\in E(\bS)$ is uniquely determined by $e\cap W$ and $e\setminus W$, where $e\setminus W\in E(\bS')$, we have $|E(\bS)|\leq2^{|W|}|E(\bS')|\leq2^{|V(\HH)|}\ex(n-|V(\HH)|,\FF)$.
\end{proof}

Given two simplicial complexes $\FF$ and $\HH$, where $\FF$ has strictly smaller dimension and $\ex(n,\FF)\ll\ex(n,\HH)$,~\Cref{obs:disjoint} implies that if $\HH$ is sought by~\Cref{p1} or~\Cref{p2}, then so is $\FF\sqcup\HH$.

In~\Cref{LB} we have actually obtained the lower bound $\ex(n,\FF) \geq \ex_s^{\cl+}(n,\FF^{s})+\sum_{r=0}^{s-1}\binom{n}{r}$ for every $s\in\{2,\dots,k\}$. We now present a construction showing that, for every $2\leq s\leq k$, there exists a $(k-1)$-dimensional simplicial complex attaining the lower bound.

\begin{obs}
\label{obs:s_bound}
Let $k\geq2$ and $s\in\{2,\dots,k\}$. For the $(k-1)$-dimensional simplicial complex $\FF=\cD(K^k_k)\sqcup\cD(K^s_{k+1})$, we have $\ex(n,\FF)=\ex_s^{\cl+}(n,\FF^{s})+\sum_{r=0}^{s-1}\binom{n}{r}$ for sufficiently large $n$.
\end{obs}
\begin{proof}[Proof of~\Cref{obs:s_bound}]
By~\Cref{new:trivial}~\eqref{(i)} we have $\ex(n,\cD(K^k_k))=\Theta(n^{k-1})$. Moreover, since $\cD(K^s_{k+1})$ contains a copy of $K_{k+1}$ in its second layer, by~\Cref{LB} and~\Cref{Zykov} we have $\ex(n,\cD(K^s_{k+1}))\geq\ex^{\cl+}_2(n,K_{k+1})=\Omega(n^k)$. Then, as $\ex(n,\cD(K^k_k))\ll\ex(n,\cD(K^s_{k+1}))$,~\Cref{obs:disjoint} implies that $\ex(n,\FF)=\ex(n,\cD(K^s_{k+1}))$. Since $\cD(K^s_{k+1})$ is trivial by~\Cref{new:trivial}~\eqref{(i)}, we have $\ex(n,\FF)=\ex^{\cl+}_s(n,K^s_{k+1})+\sum_{r=0}^{s-1}\binom{n}{r}\leq\ex_s^{\cl+}(n,\FF^{s})+\sum_{r=0}^{s-1}\binom{n}{r}$.
\end{proof}

\subsection{Connection to Berge hypergraphs}
\label{berge}
Given hypergraphs $F$ and $F'$, we say that $F'$ is a \emph{Berge copy} of $F$ (\emph{Berge-$F$} in short) if there is a bijection $f:E(F)\to E(F')$ such that for each edge $e\in E(F)$ we have that $f(e)$ contains $e$ as a subset. In other words, we can obtain $F'$ by enlarging the edges of $F$. 
The name originates from the definition of hypergraph cycles due to Berge. This was extended to arbitrary graphs by Gerbner and Palmer~\cite{gp1}, and multiple papers noted the possibility to extend it to hypergraphs. Extremal problems regarding Berge copies of hypergraphs were studied in~\cite{BGKKP}.%\vspace{1mm}

Since $F$ is a Berge copy of itself, any simplicial complex $\HH$ containing a copy of $\cD(F)$ also contains a Berge-$F$ hypergraph. Conversely, if a simplicial complex $\HH$ contains a Berge copy $F'$ of $F$, then it contains $\cD(F')$, which in turn contains a copy of $\cD(F)$. Therefore, in simplicial complexes, forbidding Berge copies of $F$ or the simplicial complex generated by $F$ is the same.%\vspace{1mm}

Previous studies focused on finding an extremal non-uniform Berge-$F$-free hypergraph that is not necessarily downward closed. However, it is immediate from the definition that if a Berge-$F$-free hypergraph contains $e$ as an edge and does not contain $e'$ as an edge, where $e'$ is a subset of $e$, then by replacing $e$ with $e'$, we obtain another Berge-$F$-free hypergraph of the same number of edges. Therefore, we can assume that the extremal hypergraph is downward closed and we have that $\ex(n,\cD(F))$ is equal to the largest number of edges in a non-uniform Berge-$F$-free hypergraph on $n$ vertices.%\vspace{1mm}

The most studied case for Berge-$F$-free hypergraphs is when $F$ is a graph. We have established in~\Cref{new:trivial}~\eqref{(i)} that for this case $\ex(n,\cD(F))=\ex_2^{\cl}(n,F)$. From the perspective of Berge copies of $F$, this result was proved in~\cite{gp2}. We remark that non-uniform Turán problems for Berge copies have received considerable attention, but most of the research focuses on weighted versions, where larger edges have larger weight.

%%%%%%%%%%%%%%%%%%%%%%%%%%%%%%%%%%%%%%%%%%%%%%%%
\subsection{Simplicial complexes on two triples}
%%%%%%%%%%%%%%%%%%%%%%%%%%%%%%%%%%%%%%%%%%%%%%%
One of the simplicial complexes addressed in~\cite{CPS} is the one generated by two disjoint triples connected by a $2$-edge. It seems natural to consider, as a warm-up problem, all simplicial complexes $\FF$ with $\cE(\FF)$ consisting of two disjoint triples and a subgraph $G$ of the complete bipartite graph $K_{3,3}$ lying cross them. Let $M^3_2+G$ denote such a simplicial complex. The result from~\cite{CPS} is then about $\ex(n,M^3_2+K_2)$. Some of our results also encompass several simplicial complexes of this type.~\Cref{new:asymptrivial} yields that $\ex(n,M^3_2+2K_2)=n^2+o(n^2)$.~\Cref{thm:non-trivial} covers $M^3_2+K_{1,3}$, $M^3_2+C_4$, $M^3_2+C_6$ and demonstrates in particular that these simplicial complexes are non-trivial. Below are our partial results on some of the remaining cases. Recall that $P_t$ denotes a graph path on $t$ edges. Recall that $G\sqcup H$ denotes the vertex-disjoint union of two hypergraphs $G$ and $H$. In particular, we write $2G$ for $G\sqcup G$.

\begin{obs}
\label{obs:short}
$\ex(n,M^3_2+P_3),\ex(n,M^3_2+2P_2),\ex(n,M^3_2+(P_2\sqcup P_1))$ are all of order $\Theta(n^2)$.
\end{obs}
\begin{proof}[Proof of~\Cref{obs:short}]
The lower bound follows by considering a simplicial complex consisting of all edges of size at most $2$ and all $3$-edges containing a fixed vertex. To see the upper bounds, observe that $\cD(TP^3_4)$ contains a copy of $M^3_2+P_3$, where $TP^3_4$ denotes a $3$-uniform tight path on $4$ edges. As a tight path is edge-degenerate, it holds by~\Cref{new:trivial}~\eqref{(i)} that $\ex(n,\cD(TP^3_4))=\Theta(n^2)$. Hence, $\ex(n,M^3_2+P_3)\leq\ex(n,\cD(TP^3_4))=\Theta(n^2)$.%\vspace{1mm}

Since $M^3_2+(P_2\sqcup P_1)$ is contained in $M^3_2+2P_2$, it suffices to show that $\ex(n,M^3_2+2P_2)= O(n^2)$. Let $F$ be the unique $3$-uniform hypergraph on $6$ vertices and $4$ edges $e_1,e_2,e_3,e_4$ such that $e_1\cap e_2=e_3\cap e_4=\emptyset$ and $e_1\cup e_2=e_3\cup e_4$. It was proved in~\cite{F} that $\ex_3(n,F)= O(n^2)$. Since $\cD(F)$ contains a copy of $M^3_2+2P_2$, for any $n$-vertex $(M^3_2+2P_2)$-free simplicial complex $\HH$, we have $|E(\HH^3)|\leq \ex_3(n,F) = O(n^2)$. One can check that if two $5$-edges intersect in at least two vertices, then their downward closure contains a copy of $M^3_2+2P_2$, so $|E(\HH^5)|=O(n^2)$. Similarly, the downward closure of two $4$-edges intersecting in exactly two vertices would contain a copy of $M^3_2+2P_2$, so $|E(\HH^4)|=O(n^2)$. Therefore, $\ex(n,M^3_2+2P_2)= O(n^2)$. 
\end{proof}

\begin{obs}
\label{obs:long}
$\ex(n,M^3_2+P_4)$ and $\ex(n,M^3_2+P_5)$ are both at most $O(n^{11/4})$. 
\end{obs}
\begin{proof}[Proof of~\Cref{obs:long}]
Let $j\in\{4,5\}$, and let $\HH$ be an $n$-vertex $(M^3_2+P_j)$-free simplicial complex. Let $m_r$ denote the number of $r$-edges in $\HH$. Since $\HH$ contains no edges of size greater than $5$ and the number of edges of size smaller than $3$ is $O(n^2)$, it suffices to show that $m_3+m_4+m_5=O(n^{11/4})$. The proof closely follows that of~\Cref{thm:non-trivial} for $\ex(n,\FF_4)$. Recall that $\FF_4$ is a copy of $M^3_2+C_6$, which contains a copy of $M^3_2+P_j$. Thus, the proof for $\ex(n,\FF_4)$ already yields that $m_3=O(n^{11/4})$ and $m_5\leq m_4$. To bound $m_4$, we apply the argument for $\ex(n,\FF_4)$, observing that the auxiliary graph $G_{x,y}$ is not only $C_4$-free but also $P_3$-free. This additional property yields that $m_4=O(n^{11/4})$.
\end{proof}

%%%%%%%%%%%%%%%%%%%%%%%%%%%%%%%%%%%%%%%%%%%%%%%%
\subsection{Possible exponent of $\ex(n,\FF)$}
%%%%%%%%%%%%%%%%%%%%%%%%%%%%%%%%%%%%%%%%%%%%%%%
Let $\FF$ be a $(k-1)$-dimensional simplicial complex on $q$ vertices. Since an $\FF$-free simplicial complex contains no edge of size $q$ or greater, we see that $\ex(n,\FF)= O(n^{q-1})$ and this bound is attained for example by $\FF=\cD(K^k_q)$. On the other hand, from the trivial lower bound we have that $\ex(n,\FF)=\Omega(n^{k-1})$.%\vspace{1mm}

In addition, for any integer $t$ with $k-1\leq t\leq q-k-1$, there is a $(k-1)$-dimensional simplicial complex $\FF$ on $q$ vertices with $\ex(n,\FF)= \Theta(n^{t})$. This can be seen by taking $\cE(\FF)$ to be a disjoint union of $K_{t+1}$ and a single $k$-edge and $q-k-t-1$ singletons. It is unclear whether there is such an $\FF$ when $q-k\leq t<q-1$.%\vspace{1mm}

A classical question in hypergraph Tur\'an theory asks which real number $t\geq0$ can occur as the exponent in the extremal number of some $k$-uniform hypergraph. An analogue of this question can be formulated for simplicial complexes.

\begin{quest}
Fix an integer $k\geq2$. For which real number $t\geq k-1$ does there exist a $(k-1)$-dimensional simplicial complex $\FF$ such that $\ex(n,\FF)=\Theta(n^t)$?
\end{quest}

The above construction, where we let $\cE(\FF)$ be a disjoint union of $K_{t+1}$ and a single $k$-edge, shows that every integer $t\geq k-1$ can occur as the exponent of $\ex(n,\FF)$ for some $(k-1)$-dimensional $\FF$. The following observation shows that the non-integer $k-1/2$ can also appear as an exponent. Given a graph $G$, the $k$-uniform \emph{suspension hypergraph} $\cS^kG$ is a hypergraph with vertex set $S\cup V(G)$ and edge set $\{S\cup e:\,e\in E(G)\}$, where $S$ is a $(k-2)$-element set disjoint from $V(G)$.

\begin{obs}
\label{exponent}
For any integer $k\geq2$, $\ex(n,\cD(\cS^kC_4))=\Theta(n^{k-1/2})$.
\end{obs}
\begin{proof}[Proof of~\Cref{exponent}]
By~\Cref{new:trivial}~\eqref{(i)} we have $\ex(n,\cD(S^kC_4))=\ex_k^{\cl+}(n,\cS^kC_4)+\Theta(n^{k-1})$. Moreover, it holds that $\ex_k^{\cl+}(n,\cS^kC_4)=\Theta(\ex(n,\cS^kC_4))$, since in any $\cS^kC_4$-free $k$-uniform hypergraph there exists no clique of order larger than $k+1$ and every edge is contained in at most one clique of order $k+1$, otherwise there would be a copy of $\cS^kC_4$. Therefore, the claim follows from a result of Mubayi~\cite[Theorem~3.1]{M}, which states that $\ex_k(n,\cS^kC_4)=\Theta(n^{k-1/2})$.
\end{proof}

Further non-integral exponents are likely to arise from considering $\cD(\cS^kG)$ for other graphs $G$. It remains unclear, however, whether a non-integral exponent greater than $k$ can occur.

\section*{Acknowledgements}
The authors would like to thank Cory Palmer, Casey Tompkins, and the Erd\H{o}s Center in Budapest for organising the ``Workshop on Generalized and Planar Tur\'an Problems'' in July 2024, where this project was initiated. They also thank Andrzej Grzesik and Daniel Johnston for fruitful discussions during the workshop. The authors acknowledge the use of AI assistance in discussing the proof of~\Cref{thm:jump} and in finding the construction recorded in~\Cref{obs:s_bound}.

\end{document}